\documentclass[twoside]{article}
\usepackage{amsmath,amssymb,latexsym}

\usepackage{amsthm}

\swapnumbers
\theoremstyle{plain}
\newtheorem{theorem}{Theorem}[section]
\newtheorem{lemma}[theorem]{Lemma}
\newtheorem{corollary}[theorem]{Corollary}
\newtheorem{proposition}[theorem]{Proposition}

\theoremstyle{definition}
\newtheorem{definition}[theorem]{Definition}

\theoremstyle{remark}

\newcommand{\reals}{\mathbb{R}}

\newcommand{\naturals}{\mathbb{N}}
\newcommand{\integers}{\mathbb{Z}}

\newcommand{\boundary}[1]{\partial#1}

\newcommand{\abs}[1]{\left\lvert#1\right\rvert} 
\newcommand{\norm}[1]{\left\lVert#1\right\rVert}

\newcommand{\tensor}{\otimes}
\newcommand{\into}{\hookrightarrow}

\DeclareMathOperator{\vol}{vol}    

\newcommand{\forget}[1]{}

{\catcode`@=11\global\let\c@equation=\c@theorem}

\begin{document}
\title{Manifolds with boundary and of bounded geometry}
\author{
Thomas Schick\thanks{
e-mail: thomas.schick@math.uni-muenster.de
\protect\\
www:~http://wwwmath.uni-muenster.de/u/lueck/
\protect\\
Fax: ++49 -251/83 38370\protect\\
Research partly funded by DAAD (German Academic Exchange Agency)
}\\
 FB Mathematik --- Universit\"at M\"unster\\
 Einsteinstr.~62 --- 48159 M\"unster\\
 Germany
}

\date{Last edited: Jan 19, 2000 --- Last complied: \today}        
\maketitle

\begin{abstract}
For non-compact manifolds with boundary  we prove that
bounded geometry defined by coordinate-free curvature bounds is equivalent to
bounded geometry defined using bounds on the metric tensor in geodesic coordinates.

We produce a nice atlas with subordinate partition of unity on
manifolds with boundary of bounded geometry, and we study the change of
geodesic coordinate maps.
\end{abstract}  

\section{Introduction}
\label{sec:intro}

Manifolds of bounded geometry arise naturally when one deals with
non-compact Riemannian manifolds, and are studied extensively in the
literature. So far, the focus was on manifolds without boundary. 

One main source of examples are coverings of compact manifolds, which
are particularly important in the context of $L^2$-cohomology
and other $L^2$-invariants. These invariants are studied frequently also for
manifolds with boundary. Therefore, it is natural to look at more
general 
manifolds with boundary and bounded geometry.

There are mainly two ways to define manifolds of bounded geometry:
either one uses bounds on the curvature (and its covariant
derivatives) ---this is the coordinate-free description--- or one uses
geodesic charts and bounds on the metric tensor and its derivatives in these
coordinates ---the coordinate approach.
A proof of the equivalence of these two definitions for manifolds without boundary
can be found in Eichhorn \cite{Eichhorn(1991a)}, using Jacobi
fields. Related but different results are obtained in \cite{Atiyah-Bott-Patodi(1973)}  using
synchronous frames. The case of
manifolds with boundary causes additional technical difficulties and
seems not to be covered in
the literature. Therefore we give a proof here, using synchronous frames.

Dealing with manifolds with boundary, in addition to the usual
requirements in the interior we must impose boundary regularity
conditions. These involve the second fundamental form (in the coordinate-free
description) or special charts for the boundary (in the coordinate
description).

In the last section, we show that the functions given by a change of geodesic coordinates
 and their derivatives admit uniform bounds on manifolds of
bounded geometry. And we provide one technical tool, namely a nice atlas with
subordinate nice partition of unity. This was
introduced and used by Shubin \cite[1.2 and 1.3]{Shubin(1992b)} if the
boundary is empty.

This paper grew out of part of the Dissertation \cite{Schick(1996)} of the author, and the results obtained here  are used in
\cite{Schick(1998d)}. I 
thank my advisor, Prof.~Wolfgang L\"uck, for his constant support and
encouragement. I also thank the referees for valuable comments and
suggestions concerning the exposition of the paper.

\section{Coordinate-free versus coordinate-wise curvature bounds}
\label{sec:curv}

\begin{definition}\label{basicnotation}
  On a Riemannian manifold $(M^m,g)$ with boundary $\boundary M$, $R$
  denotes the curvature
  tensor of $M$, $l$ the second fundamental form of $\boundary M$, and
  $\bar R$ the curvature tensor of $\boundary M$ (with its induced
  metric). The (Levi-Civita)-covariant derivative of $M$ is denoted
  with $\nabla$, the one of $\boundary M$ with $\bar\nabla$. We
  use $\nu$ for the unit inward normal vector field at $\boundary M$.

  If not stated otherwise, a manifold $M$ will always have dimension
  $m$.

  Given an open subset $U\subset M$ and  a chart
  $x=(x_1,\dots,x_m):U\to\reals^m$, we consider the
  corresponding derivations $\frac{\partial}{\partial x_i}$ as
  derivations on $U$, or as elements in the tangent
  bundle $TM$. We abbreviate
  $\partial_i:=\frac{\partial}{\partial x_i}$.
  We let $g_{ij}:=g(\partial_i,\partial_j)$ be the metric tensor in
  the given coordinates and $g^{ij}$ be the coefficients of the
  inverse matrix.

  We use the notation of multi-indices throughout:
  Let $\alpha=(\alpha_1,\dots,\alpha_m),
  \beta=(\beta_1,\dots,\beta_m)$ be multi-indices (with
  $\alpha_i,\beta_i\in\naturals\cup\{0\}$). Then 
  \begin{equation*}
   D^\alpha:= D_x^\alpha:= \frac{\partial^{\alpha_1}}{\partial
    x_1^{\alpha_1}}\dots\frac{\partial^{\alpha_m}}{\partial
    x_m^{\alpha_m}} ;
  \end{equation*}
  and we set $\beta\le\alpha$ if and only if $\beta_i\le \alpha_i$ for
  $i=1,\dots,n$. Define $\abs{\alpha}:=\sum_{i=1}^m\alpha_i$.

  For $V\subset M$ and $r>0$ set $U_r(V):=\{x\in M|\; d(x,V)<r \}$. For
  $p\in M$ set $B(p,r):=U_r(\{p\})$. If
  $p\in\boundary M$, $(B(p,r)\subset\boundary M)$ means the
  corresponding set for $\boundary M$ with the induced Riemannian
  metric.

  We use the normal geodesic flow $K:\boundary M\times [0,\infty)\to
  M: (x',t)\mapsto \exp^M_{x'}(t\nu_{x'})$. For $p\in\boundary M$ set
  $Z(p,r_1,r_2):= K( (B(p,r_1)\subset \boundary M)\times
  [0,r_2))\subset M$.

  Set $N(s):= K(\boundary M\times [0,s])$ if $s\ge 0$.
\end{definition}

\begin{definition}\label{def_bounded} Suppose $M$ is a manifold with
  boundary $\boundary M$ (possibly empty). It is of {\em
  (coordinate-free defined) bounded geometry} if the following holds:
\begin{itemize}
\item[(N)] {\em Normal collar:} there exists $r_C>0$ so that the
  \emph{geodesic collar} 
\begin{equation*}  \boundary M \times [0,r_C) \to M: (t,x)\mapsto \exp_x(t\nu_x) \end{equation*}
is a diffeomorphism onto its image ($\nu_x$ is the unit inward normal
vector).
\item [(TIC)] The {\em injectivity radius $r_{inj}(\boundary M)$ of $\boundary M$} is positive.
\item[(I)] {\em Injectivity radius of $M$:} There is $r_i>0$ so that
  if $r\le r_i$ then
  for $x\in M-N(r)$ the exponential map is a diffeomorphism on
  $B(0,r)\subset T_xM$. Hence, if we identify $T_xM$ with $\reals^m$ via
  an orthonormal frame we have \emph{Gaussian coordinates} $\reals^m\supset B(0,r)\stackrel{\exp^M_x}{\to} M$ around every point in $M-N(r)$.
\item[(B)] {\em Curvature bounds:} For every $k\in\naturals$ there is $C_k>0$
so that $\abs{\nabla^i R}\le C_k$ and $\abs{\bar\nabla ^i l}\le C_k$ for
$ 0\le i\le k$.
\end{itemize}
\end{definition}

The injectivity radius and curvature bounds are what one is used to
for manifolds without boundary (compare e.g.~\cite[Section
3]{Cheeger-Gromov-Taylor(1982)}). The embedding of the boundary
 is described by the second fundamental form. Because the injectivity radius does not make sense near the boundary, we  replace it by the geodesic collar.

To give the coordinate-wise definition of bounded geometry, we have to explain
which charts we want to use:
\begin{definition} Let $M$ be a Riemannian manifold with boundary
  $\boundary M$.
Fix $x'\in\boundary M$ and
an orthonormal basis of $T_{x'}\boundary M$ to identify
$T_{x'}\boundary M$ with $\reals^{m-1}$. For $r_1,r_2>0$
sufficiently small (such that the following map is injective) define
{\em normal collar coordinates}
\begin{equation*} \kappa_{x'}: \underbrace{B(0, r_1)}_{\subset
    \reals^{m-1}}\times [0,r_2)\to M: (v,t)\mapsto
  \exp^M_{\exp_{x'}^{\boundary M}(v)}(t\nu) .\end{equation*}
(We compose  the exponential maps of $\boundary M$ and  of $M$, and
  $\nu$ is the inward unit normal vector 
field). The tuple $(r_1,r_2)$ is called the \emph{width of the normal
  collar chart} $\kappa_{x'}$.

We adopt the convention that the boundary defining coordinate is the
last (i.e.~$m^{\text{th}}$) coordinate.

For $x\in M-\boundary M$ and $r_3>0$ sufficiently small the
exponential map yields {\em Gaussian coordinates} (identifying $T_xM$
with $\reals^m$ via an orthonormal base)
\begin{equation*} \kappa_x: B(0,r_3)\to M: v\mapsto \exp_x^M(v)
  . \end{equation*}
We call $r_3$ the \emph{radius of the Gaussian chart} $\kappa_x$.

We use the common name {\em normal coordinates} for
normal collar coordinates as well as Gaussian coordinates.

\end{definition}

\begin{definition} \label{coorddef}
A Riemannian manifold $M$ with boundary $\boundary M$ has {\em
    (coordinate-wise defined) bounded
    geometry} if and only if (N), (IC), (I) of Definition \ref{def_bounded} hold and (instead of (B))\\
\begin{itemize}
\item[(B1)] There exist $0<R_1\le r_{inj}(\boundary M)$, $0<R_2\le r_C$
  and $0<R_3\le r_i$ and constants $C_K>0$ (for each $K\in\naturals$)
  such that whenever we have normal boundary
  coordinates of width $(r_1,r_2)$ with $r_1\le R_1$ and $r_2\le R_2$,
  or Gaussian coordinates of radius $r_3\le R_3$ then
   in these coordinates
\begin{equation*} \abs{D^\alpha g_{ij}}\le C_K\quad\text{and}\quad
\abs{D^\alpha g^{ij}}\le C_K\qquad\forall\abs{\alpha}\le K . \end{equation*}
\end{itemize}

The numbers $R_1$, $R_2$, $R_3$ and $C_K$ are called the bounded
geometry constants of $M$.
\end{definition}

The main result of the paper is the following:
\begin{theorem}\label{main}
 Let $(M^m,g)$ be a Riemannian manifold with boundary $\boundary M$.
To given $C>0$, $k\in\naturals$, and dimension $m$ there are
$R_1,R_2,R_3>0$ and $D>0$
such that the
following holds:
\begin{itemize}
\item[(a1)] If $x\in M- \boundary M$, $0<r_3\le R_3$ and
  $\kappa_x:B(0,r_3)\to (M-\boundary M)$ is a Gaussian chart, and if
  $\abs{\nabla^i R}\le C$ for $i=0,\dots,k$ on the image of $\kappa_x$
  then in these coordinates
  \begin{equation*}
    \abs{D^\alpha g_{ij}}\le D \quad\text{and}\quad
\abs{D^\alpha g^{ij}}\le D\qquad\text{whenever }\abs{\alpha}\le k.
  \end{equation*}
\item[(a2)]
  If on the other hand 
\begin{equation*}
    \abs{D^\alpha g_{ij}}\le C\quad\text{and}\quad
\abs{D^\alpha g^{ij}}\le C\qquad\text{for }\abs{\alpha}\le k+2
\end{equation*}
 then on the image of $\kappa_x$ we have
 \begin{equation*}
\abs{\nabla^i R}\le D\qquad\text{for } i=0,\dots,k .
\end{equation*}
\item[(b1)] If $x'\in\boundary M$, $0<r_1\le R_1$, $0<r_2\le R_2$ and
  $\kappa_{x'}: B(0,r_1)\times[0,r_2)\to M$ is a normal boundary
  chart, and if $\abs{\nabla^i R}\le C$ and $\abs{\bar\nabla^i l}\le
  C$ for $i=0,\dots,k$ on the image of $\kappa_{x'}$, then in these
  coordinates we get
  \begin{equation*}
    \abs{D^\alpha g_{ij}}\le D\quad\text{and}\quad
\abs{D^\alpha g^{ij}}\le D\qquad\text{whenever }\abs{\alpha}\le k .
  \end{equation*}
\item[(b2)] If, on the other hand,
\begin{equation*}
    \abs{D^\alpha g_{ij}}\le C\quad\text{and}\quad
\abs{D^\alpha g^{ij}}\le C\qquad\text{for }\abs{\alpha}\le k+2
\end{equation*}
 then, on the image of
 $\kappa_{x'}$, 
 \begin{equation*}
\abs{\nabla^i R}\le D\text{ and }\abs{\bar\nabla^i l}\le
 D\qquad\text{for }
 i=0,\dots,k.
\end{equation*}
\item[(c)] $M$ has (coordinate-wise defined) bounded geometry if and only if it has (coordinate-free defined)
bounded geometry. In particular, we can drop the prefix in
notation. The bounded geometry constants of Definition
\ref{coorddef} can be chosen to depend only on $r_i$, $r_C$,
$r_{inj}(\boundary M)$ and $C_k$ of Definition \ref{def_bounded}.
\end{itemize}
\end{theorem}
Observe that (c) follows from
(a1)-(b2). Moreover, (a2) and
(b2) are immediate consequences of the formulas for $R$ and
$l$ (and their covariant derivatives) in local coordinates in terms of
$g_{ij}$, $g^{ij}$ and their
partial derivatives (compare 2.54, 3.16 and 5.1 of
\cite{Gallot-Hulin-Lafontaine(1987)} --- note that our charts near the
boundary are adapted to the embedding $\boundary M\into M$). The
statement (a1) about internal points is already included in
\cite[Theorem A and Proposition 2.3]{Eichhorn(1991a)}. It remains to
establish (b1). Since in the course of this proof we have to set up
most of the notation necessary for the synchronous-frame-proof of
(a1), we include a
complete proof also of (a1).

The proof is done in four steps. First, we give the argument for
$k=0$, using the Rauch
comparison theorem. Secondly, we prove (a1). In the third step, we
establish bounds on the curvature tensor of the boundary. Last, we
derive (b1).

\medskip
\noindent\textbf{Step 1: Proof of Theorem \ref{main}(a1) and
  \ref{main}(b1) for $k=0$}
\begin{proposition}\label{rauch}
Suppose we are in the situation of Theorem \ref{main}(a1) or (b1) and
$k=0$. Suppose $(x_i):=\kappa^{-1}:U\subset M \to \reals^m$ is the
normal coordinate system. There are $R_1,R_2,R_3>0$ and $C_1,C_2>0$
(depending only on
$C$ and $m$) such that if for width or radius we have
$(r_1,r_2)\le (R_1,R_2)$ or $r_3\le R_3$, respectively, then
\begin{equation}\label{b_Rauch}
 C_1\le \abs{\sum\lambda_i\frac{\partial}{\partial x_i}}_{TM}\le
 C_2,\qquad\text{if }\sum\lambda_i^2=1,
\end{equation}
where
$\abs{v}_{TM}:=\sqrt{g(v,v)}$ for $v\in TM$.

Moreover, $g_{ij}$ and $g^{ij}$ are bounded with a bound depending only
on $C$ and the dimension.
\end{proposition}
The numbers $R_1$, $R_2$ and $R_3$ of Theorem \ref{main} are
determined by Proposition \ref{rauch} and (IC), (I), (N).

\begin{proof}
The last statement is a reformulation of Inequality
\eqref{b_Rauch}. To prove \eqref{b_Rauch}, we apply Warner's generalization
of the Rauch  comparison theorem \cite[4.3]{Warner(1966)}. We compare
with two complete manifolds of constant sectional curvature $-C$ and
$C$, respectively. To compare with normal collar coordinates, choose a
hypersurface in this manifold so that all the eigenvalues of its
second fundamental form at one (comparison) point are equal to $C$ in
the first case and to $-C$ in the second case. Inequality
\eqref{b_Rauch} for vectors orthogonal to ${\mathcal{R}}:=\sum
x_i\partial_i$ (in Gaussian coordinates), or orthogonal to
$\partial_m$ (in normal boundary coordinates) is just the statement of
the comparison theorem, with
$C_1$ and $C_2$ depending only on the manifold we compare
with (i.e.~on $C$ and on $m$). Here, $r_1$, $r_2$ and $r_3$ must be sufficiently small
(again depending only on the manifolds we compare with). The
comparison theorem says nothing about ${\mathcal{R}}$ or about
$\partial_m$, respectively. But for these vectors Euclidean length
and length in $TM$ as well as the orthogonal complements coincide by
the following Proposition \ref{Gauss_lemma}. Therefore, the
inequality is true in general.
\end{proof}

In the proof of Proposition \ref{rauch} we used  the Gauss lemma:
\begin{proposition}\label{Gauss_lemma}
Let $(M,g)$ be a Riemannian manifold and $\exp:B(0,R)\to M$ a Gaussian chart.
Pull the metric $g$ back to $B(0,R)$. Then
$g({\mathcal{R}},{\mathcal{R}})=r^2$, (${\mathcal{R}}=\sum_i
x_i\partial_i$), and $g({\mathcal{R}},v)=0$ if and only if $v$
is a tangent vector to a sphere with center the origin $0$.\\
Let $K:\boundary M\times[0,r_C)\to M$ be the geodesic collar and pull
$g$ back to $\boundary M\times [0,r_C)$. Then
$g(\partial_m,\partial_m)=1$ and $g(\partial_m,v)=0$ if and only if
$v$ is tangent to a translate $\boundary M\times
\{t\}$.
\end{proposition}
\begin{proof}
Compare \cite[2.93]{Gallot-Hulin-Lafontaine(1987)} --- the proof there
works also for the collar.
\end{proof}

\medskip
\noindent\textbf{Step 2: Proof of \protect\ref{main}(a1).}

Suppose we are in the situation of \ref{main}(a1) with $p\in
M-\boundary M$ and Gaussian coordinates
$x=(x_1,\dots,x_m)=\kappa_p^{-1}:B(p,r_3)\to\reals^m$.
We will state a (differential) equation for $g_{ij}$ in terms of the
curvature tensor, so that a bound on partial derivatives of the
components of the curvature tensor will give corresponding bounds for
the metric. Partial  and covariant derivative are related by the
Christoffel symbols, so we will compute them, too.

Choose an orthonormal base
$\{s_i\}$ for $T_pM$. Using parallel transport along geodesics
emanating from $p$, construct a synchronous orthonormal frame
$\{s_i(x)\}$ of the tangent space restricted to $B(p,r_3)$. Let
$\{\theta^i\}$ be the frame of 1-forms dual to
$\{s_i\}$ (therefore orthonormal). The connection forms
$\theta_j^i$ for this frame are defined by 
\begin{equation*}\nabla s_j=\sum_i\theta^i_js_i,
\end{equation*}
 with associated Christoffel symbols $\Gamma^i_{jk}$ and curvature tensor $R^i_{jkl}$ given by
 \begin{equation*}  \theta^i_j =\sum_k \Gamma^i_{jk}
dx_k;\qquad
  d\theta^i_j-\sum_k\theta^i_k\wedge\theta^k_j 
=\sum_{k,l}R^i_{jkl}dx_k\wedge dx_l. 
\end{equation*}
We can express the curvature entirely in terms of $s_i$ and
  $\theta^i$, which defines $K^i_{jkl}$:
\begin{equation*} d\theta^i_j-\sum_k\theta^i_k\wedge\theta^k_j = \sum_{k,l}K^i_{jkl}\theta^k\wedge\theta^l. \end{equation*}
Define functions $a^i_{j}$ and $b^i_{j}$ via the equations
\begin{equation}\label{defa}
 \theta^i=\sum_j a^i_j dx_j; \qquad dx_i=\sum_j b^i_j
\theta^j. \end{equation}
\begin{equation}\label{R_K}
  \begin{split}
\text{Then}\quad R^i_{jkl} &=\sum_{\alpha,\beta}
 K^i_{j\alpha\beta}a^\alpha_k a^\beta_l\quad\text{and}\quad 
 g_{ij}=\sum_\alpha a^\alpha_i a^\alpha_j;\quad g^{ij}=\sum_\alpha
b^i_\alpha b^j_\alpha .
\end{split}\end{equation}
As matrix, $(g_{ij})$ is the product of $(a^i_j)$
and its adjoint, and accordingly for $(g^{ij})$ and $(b^i_j)$. Hence
\begin{lemma}\label{abounds}
  There are bounds on $a^i_j$ and $b^i_j$ corresponding to the bounds
  on $g_{ij}$ and $g^{ij}$ given by Proposition \ref{rauch}.
\end{lemma}

The Christoffel symbols
$\tilde\Gamma^i_{jk}$ of the covariant differentials of
$\partial_i$ are given by
\begin{equation*} \nabla_{\partial_k}\partial_j = \sum_i
\tilde\Gamma^i_{jk}s_i. \end{equation*}
Dualizing \eqref{defa} we see that $\partial_j=\sum_\alpha a^\alpha_j
s_\alpha$, hence
\begin{equation}\label{tgamma_gamma}
 \tilde\Gamma^i_{jk} = \partial_k a^i_j + \sum_\alpha
 a^\alpha_j\Gamma^i_{\alpha k}.
\end{equation}

Atiyah, Bott, and Patodi \cite[a6 and
a10]{Atiyah-Bott-Patodi(1973)} derive the following equations (note
that our definition of $R^i_{jkl}$ takes care of the
problems described in \cite{Atiyah-Bott-Patodi(1975)}), where
${\mathcal{R}}=\sum_i
x_i\partial_i$:
\begin{gather}\label{Gamma}
{\mathcal{R}}\Gamma^i_{jk}+\Gamma^i_{jk}=\sum_l
2x_lR^i_{jkl}\qquad\forall i,j,k ; \\
\label{gij}
({\mathcal{R}}^2+{\mathcal{R}})a^i_l=-2\sum_{j,k}R^i_{jkl}x_j
x_k\qquad\forall i,l .
\end{gather}
Set $f_x(t):=t\Gamma^i_{jk}(tx)$. Let $'$ denote differentiation with
respect to $t$. Then
\begin{gather} f_x'(t)=\Gamma^i_{jk}(tx)+t\sum_l
x_l\partial_l\Gamma^i_{jk}(tx) \stackrel{\eqref{Gamma}}{=}\sum_l
2t\cdot x_lR^i_{jkl}(tx). \nonumber \\
  \implies \quad \Gamma^i_{jk}(x)=f_x(1)=\int_0^1\sum_l2\tau
  x_lR^i_{jkl}(\tau x)\;d\tau \qquad\text{and} \nonumber \\
  \label{bgamma}
D^\alpha_x\Gamma^i_{jk}(x)=\int_0^1\tau^{\abs{\alpha}}\left(D^\alpha_x(x\mapsto
\sum_l x_lR^i_{jkl}(x))\right)(\tau x)\;d\tau .
\end{gather}

Set $f_{il}(t,x):=a^i_l(t x)$. 
Then $t {f_{ij}}'(t,x)={\mathcal{R}} a^i_j(tx)$ and
$t^2f_{il}''(t,x)+tf_{il}'(t,x)={\mathcal{R}}^2a^i_l(tx)$. By \eqref{gij}
\begin{equation*} t^2f_{il}''+2tf_{il}' =
-2t^2\sum_{k,j}R^i_{jkl}(tx)x_j x_k . \end{equation*}
With $w_{il}(t,x):=t^2f_{il}'(t,x)$ we get
$w_{il}'=t^2f_{il}''+2tf_{il}'$. Since $w_{il}(0)=0$,
\begin{equation*}t^2f_{il}'(t,x)=-2\int_0^t\tau^2\sum_{j,k}R^i_{jkl}(\tau
x)x_j x_k\;d\tau \quad\stackrel{\tau=tu}{\implies} \end{equation*}
\begin{equation}\label{K_eq}
f_{il}'(t,x)=-2t\int_0^1
u^2\sum_{j,k,\alpha,\beta} K^i_{j\alpha\beta}(tux)a^\alpha_k(tux)
a^\beta_l(tux) x_j x_k\; du .
\end{equation}

Now we are in the position to explain how the bounds on $R$ and its
covariant derivatives up to order $k$ give rise to bounds on $g_{ij}$,
$g^{ij}$ and their partial derivatives up to order $k$. Because of
\eqref{R_K} we can consider $a^i_j$ and $b^i_j$ instead of the
metric tensor. Moreover, the case $k=0$ is done by Proposition
\ref{rauch}.

\begin{lemma}\label{inv_mat}
Let $A$, $B$ be matrix valued functions which are inverse to each other. Then
\begin{equation*} \frac{\partial}{\partial x_i}B=
  \frac{\partial}{\partial x_i}(A^{-1})=-
  A^{-1}(\frac{\partial}{\partial
    x_i}A)A^{-1}=-B(\frac{\partial}{\partial x_i} A)B.
 \end{equation*}
Iterated application of this and of the product rule yields
\begin{equation*} D^\alpha_x B =P_\alpha(B, D^\beta_x A;\;\beta\le
  \alpha), 
\end{equation*}
where $P_\alpha$ is a fixed
polynomial in non-commuting variables. Bounds for the partial
derivatives of $A$ up to
order $k$ and on $B$ yield bounds for the partial
derivatives of $B$.
\end{lemma} 

 Lemma \ref{inv_mat} applies to the matrices $A=(a^i_j)$ and
 $B=(b^i_j)$. Moreover, by Proposition \ref{rauch}, we have a
 bound for $(b_{ij})$. Hence it remains to find bounds for the
 derivatives of $(a^i_j)$. 

\begin{lemma}\label{Rder} For $\alpha=(\alpha_1,\dots,\alpha_n)$ there
 is a polynomial $P_{\alpha,ijkl}$
(only depending on $\alpha,i,j,k,l$) in partial derivatives up to
order $(\abs{\alpha}-1)$ of $K^*_{***}$, $\Gamma^*_{**}$, and $a^*_*$
such that as functions on the set $B(p,r_3)$
\begin{equation*}
(\nabla_{\partial_{1}})^{\alpha_1}\dots(\nabla_{\partial_{n}})^{\alpha_n}R(s_i,s_j,\partial_k,\partial_l)=
D_x^\alpha K^i_{jkl} + P_{\alpha,ijkl} . \end{equation*}

\end{lemma}
\begin{proof}
This follows from the formula for covariant differentials in
coordinates. Note that for $\abs{\alpha}=1$ only $\Gamma^i_{jk}$ shows
up (since $K^i_{jkl}$ is defined entirely in terms of
$s_i$). But if we iterate the covariant differentials, we have to take
into account that we contracted $\nabla R$ with $\partial_i$ and not
with $s_i$. This yields (via $\nabla\partial_i$) $\tilde\Gamma^i_{jk}$
and, since we iterate the covariant differentials, their partial
derivatives up to order $\abs{\alpha}-2$. Since
$\tilde\Gamma^i_{jk}=\partial_k a^i_j +\sum_\alpha
a^\alpha_j\Gamma_{\alpha k}^i$, the result follows.
\end{proof}

Now we proceed by induction on the order of derivatives $\abs{\alpha}$.
For $\abs{\alpha}=0$ observe that by assumption we have a
bound on the
curvature. Since $\{s_i\}$ is orthonormal 
this gives bounds on $K^i_{jkl}$. By Proposition
\ref{rauch} the same is true for $a^i_j$.

Assume by induction that for $r\ge 0$ we
have found bounds on the partial derivatives up to order $r$
of $K^i_{jkl}$ and  $a^i_j$ and on the derivatives up to order $(r-1)$
of $\Gamma^i_{jk}$. The assumptions of the Theorem give bounds on
$\abs{R},\dots,\abs{\nabla^{r+1}R}$.

From equation \eqref{R_K}, relating $K^i_{jkl}$ and $R^i_{jkl}$, we
get bounds on the partial derivatives up to order $r$ of
$R^i_{jkl}$. Then Equation \eqref{bgamma} yields bounds for
the  derivatives of order $r$ of $\Gamma^i_{jk}$.  Lemma
\ref{Rder} and the bound on $\nabla^{r+1} R$ yield bounds on
$(r+1)$-order partial derivatives of $K^i_{jkl}$ (since by Proposition
\ref{rauch} the length
of $\partial_i$ is controlled).
In all instances the new bounds are given in terms of the old ones.

It
remains to deal with the derivatives of order $(r+1)$ of $a^i_j$.
Remember Equation \eqref{K_eq} for $f_{il}(t,x)=a^i_l(t x)$:
\begin{equation*}f_{il}'(t,x)=-2t\int_0^1 u^2\sum_{j,k,\alpha,\beta}
\left(K^i_{j\alpha\beta}a^\alpha_k a^\beta_l\right)(tux) x_j x_k\; du.
 \end{equation*}
Let $\alpha$ be a multi-index with $\abs{\alpha}=r+1$. We
differentiate the equation with respect to $x$ to get an equation for
$D^\alpha_x f_{il}(t,x)=t^{\abs{\alpha}}(D^\alpha a^i_l)(tx)$. This yields
\begin{multline}\label{integral}
(D^\alpha_xf_{il})'(t,x)=
-2t\int_0^1\bigl
( u^2\sum_{j,k,\beta,\gamma}K^i_{j\beta\gamma}(tux)\cdot \\
((D^\alpha_x f_{\beta k}) f_{\gamma l} +
+(D^\alpha f_{\gamma l})f_{\beta k})(tu,x) x_j x_k\;du\bigr)
-2t\int_0^1 P_\alpha\;du  .
\end{multline}
Here $P_\alpha$ is a polynomial in $t$, $u$, $x$, partial derivatives
up to order $(r+1)$ of $K^*_{***}$ at $tux$, and partial derivatives
up to order $r$ of $f_{**}$ at $(tu,x)$. The left and right hand side
of \eqref{integral} are equal as function of $x$ and $t$. The
induction hypothesis implies for $0\le
t\le 1$ with suitable $C_1,C_2>0$ the inequality
\begin{equation}\label{ader} \abs{(D^\alpha_x f_{ij})'(t,x)}\le C_1 \sup_{0\le\tau\le t}\{ \abs{D^\alpha_x f_{ij}(\tau,x)}\} +C_2 , 
\end{equation}
Moreover, $D^\alpha f(0,x)=0$ since
$\abs{\alpha}\ge 1$.

Let $h(t):=C_2(\exp(C_1t)-1)/C_1$ be the unique solution of $h'(t)=C_1
h(t)+ C_2$ with
$h(0)=0$. This is a positive monotonous increasing function, with an
explicit bound $h(t)\le C:=C_2(\exp(C_1)-1)/C_1$ for $0\le t\le 1$.\\
Abbreviate $u^i_j(t):=D^\alpha_x f_{ij}(t,x)$.
We will prove $\abs{u^i_j(t)}\le h(t)$ and therefore
\begin{equation}\label{wantit}
 \abs{D^\alpha a^i_j(x)}=\abs{u^i_j(1)}\le h(1)\le C.
\end{equation}
This then finishes the induction step.
To show $\abs{u^i_j(t)}\le h(t)$, let
$h_n$ be the unique
solution of
\begin{equation*}h_n'(t)=C_1 h_n(t) +C_2+1/n\quad\text{with
    $h_n(0)=0$}.\end{equation*}
 Then $h_n(t)\xrightarrow{n\to\infty} h(t)$ uniformly for $0\le t\le
 1$. Therefore, it suffices to show $\abs{u^i_j(t)}\le h_n(t)$.
For a contradiction, assume
$\abs{u^i_j(t)}>h_n(t)$
for some $n$ and $t$. Set $t_0:=\inf_{0\le
t}\{\abs{u^i_j(t)}>h_n(t)\}$. Then $\abs{u^i_j(t_0)}=h(t_0)$, since
$u_i^j(0)=0=h(0)$, $h_n$ is monotonous, and $\abs{u^i_j(t)}\le h_n(t)$
$\forall t\le t_0$. Consequently, $\sup_{t\le
 t_0}\abs{u^i_j(t)}=\abs{u^i_j(t_0)}$. Then \eqref{ader}
 shows
\begin{equation*} \abs{(u^i_j)'(t_0)}\le C_1\abs{u^i_j(t_0)}+C_2 <
 h_n'(t_0).  \end{equation*}
Moreover, $d/dt\abs{u^i_j(t_0)}\le \abs{(u_j^i)'(t_0)}$ (compare
\cite[III.3.2]{Hartmann(1964)} for the
difficult case $u^i_j(t_0)=0$ ---$d/dt$ is understood to be the right
derivative). It follows
$\abs{u^i_j(t)}<h_n(t)$ for $t\in [t_0,t_0+\epsilon)$ and $\epsilon>0$
sufficiently small. But this contradicts the choice of $t_0$.

\medskip
\noindent\textbf{Step 3: Curvature of $\boundary M$.}

We adopt the notation of Definition \ref{basicnotation}.

In the following we consider $(0,p)$-tensors $T$ on $M$ and their
restriction to $\boundary M$, given by the inclusion $T\boundary
M\into TM$. We will use the same notation for $T$ and its restriction, the
meaning will be clear from the context.

We compute the covariant derivatives $\bar\nabla^k\bar R$ using the following rules:
\begin{lemma}\label{rules}
Suppose $T$ is a $(0,q)$-tensor on $M$, $S$ a $(0,p)$-tensor on
$\boundary M$, and $S^{*_1}$ the $(1,p-1)$-tensor on $\boundary M$
given by $g(S_x^{*_1}(v_2,\dots,v_p),v_1)=S_x(v_1,\dots,v_p)$ for
$v_1,\dots,v_p\in T_x\boundary M$, where $x\in\boundary M$ and $S_x$,
$S_x^{{*_1}}$ are the values of $S$ and $S^{{*_1}}$, respectively, at $x$. Let $\sigma$ be a permutation (operating on a multiple tensor product by permutation of the factors) with $\sigma^{-1}(1)\le p$.\\
Let $c$ denote the contraction of a $(0,r)$-tensor with a
$(1,s)$-tensor which contracts the $r$-th entry of the $(0,r)$-tensor.
The covariant derivative is understood to be a map  $\nabla:C^\infty(E)\to
C^\infty(E\tensor T^*M)$. Then the following holds:
\begin{enumerate}
\item\label{f1} $\bar\nabla T=\nabla T -\sum_i c(T\tensor l\circ
  \sigma_i,\nu)$, where $\sigma_i$ are appropriate permutations.
\item\label{f2} $\bar\nabla((T\tensor S)\circ \sigma)=((\bar\nabla T)\tensor
  S)\circ\sigma + (T\tensor\bar\nabla S)\circ\sigma$.
\item\label{f3} $\bar\nabla c((T\tensor S)\circ\sigma,\nu) = c((\nabla
  T)\tensor
S\circ\sigma',\nu)+ c(T\tensor\bar\nabla S\circ\sigma'',\nu) + \sum_i
c(c(T\tensor S\circ\sigma,\nu)\tensor l\circ\sigma_i,\nu) + c(T\tensor
S\circ\sigma,l^{*_1})$, with $\sigma'$, $\sigma''$, and $\sigma_i$ appropriate
permutations. 
\item\label{f4} $\bar\nabla c(T,(\bar\nabla^k l)^{*_1})=c(\bar\nabla
T,(\bar\nabla^kl)^{*_1}) + c(T,(\bar\nabla^{k+1}l)^{*_1})$. 
\end{enumerate}
\end{lemma}



\begin{proof}
Formulas \ref{f2}. and \ref{f4}. are well known. Let
$v_1,\dots,v_p$ and $X$ be vector fields on $\boundary M$ For \ref{f1}.
we compute:
\begin{equation*}\begin{split}
&\bar\nabla T(v_1,\dots,v_p,X) \\
&= X.T(v_1,\dots,v_p) - T(\bar\nabla_Xv_1,\dots,v_p) -\dots-T(v_1,\dots,\bar\nabla_X v_p)\\
&\stackrel{\bar\nabla_X Y=\nabla_X Y-l(X,Y)\nu}{=} X.T(v_1,\dots,v_p) - T(\nabla_X v_1,\dotsc) - \dots\\
& \quad + T(\nu,v_2,\dotsc)l(v_1,X)+\dots+T(v_1,\dots,v_{p-1},\nu)l(v_p,X)\\
&= \nabla T(v_1,\dots,v_p,X) +\sum_i c(T\tensor
l\circ\sigma_i,\nu)(v_1,\dots,v_p,X) .
\end{split}\end{equation*}
For \ref{f3}. set $v_1:=\nu$ and calculate:
\begin{equation*}\begin{split}
  &\bar\nabla_X c(T\tensor S\circ\sigma,\nu)(v_2,\dots,v_{p+q})
  \stackrel{v_1=\nu}{=} \\
  &=(X.T(v_{\sigma1},\dots,v_{\sigma p}))S(v_{\sigma(p+1)},\dots,v_{\sigma(p+q)})\\
  &\quad + T(v_{\sigma1},\dots,v_{\sigma
    p})(X.S(v_{\sigma(p+1)},\dots,v_{\sigma(p+q)}))\\
  &\quad - \sum_{\substack{i=1\\\sigma i\ne 1}}^{p+q} T\tensor S(v_{\sigma 1},\dots,\bar\nabla_X v_{\sigma i},\dots,v_{\sigma(p+q)})\\
&= c(T\tensor\bar\nabla_X S\circ\sigma,\nu)(v_2,\dots,v_{p+q})\\
&\quad + \left(X.T(v_{\sigma 1},\dots,v_{\sigma p})- \sum_{i=1}^pT(v_{\sigma 1},\dots,\nabla_Xv_{\sigma i},\dots,v_{\sigma p})\right)S(\dotsc)\\
&\quad -\sum_{\substack{i=1\\\sigma(i)\ne 1}}^p
\underbrace{T(v_{\sigma 1},\dots,l(X,v_{\sigma i})\nu,\dots,v_{\sigma
    p})}_{=T(\cdots)l(X,v_{\sigma i})}S(\dotsc) \\
&\quad + T(v_{\sigma 1},\dots,\nabla_X\nu,\dots,v_{\sigma p})S(\dotsc)\\
&= c(T\tensor\bar\nabla_X S\circ\sigma,\nu)(\dotsc) + c((\nabla_X T)\tensor S\circ\sigma,\nu)(\dotsc) \\
&\quad -\sum_i c(c(T\tensor S\circ\sigma,\nu)\tensor
l\circ\sigma_i,\nu)(\dotsc,X) + c(T\tensor
S\circ\sigma,\nabla_X\nu)(\dotsc) .
\end{split}\end{equation*}
If $Y\in C^\infty(T\boundary M)$ then
\begin{equation*}\begin{split} & 0 = X.g(\nu,Y)=g(\nabla_X\nu,Y)+g(\nu,\nabla_X Y)\\
&\implies g(\nabla_X\nu,Y)= l(X,Y)=l(Y,X)\\
& 0 = X.g(\nu,\nu)=2g(\nabla_X\nu,\nu) .\\
\implies \quad &\nabla_X\nu = l^{*_1}(X)\implies \nabla\nu = l^{*_1}
.
\end{split}\end{equation*}
This finishes the proof.
\end{proof}

\begin{corollary}\label{boundarybounds}
$\bar\nabla^k\bar R$ is a finite sum of tensor products and possibly
iterated contractions, composed with permutations, involving
(i) $\nabla^j R$ for $j\le k$;  (ii) $\bar\nabla^j l$ for $  j<k$;
    (iii) $(\bar\nabla^j l)^{*_1}$ for $j< k-1$; and   (iv) $\nu$.\\
Bounds for the building blocks (i) and
(ii) yield a  bound for $\bar\nabla^k\bar R$.
\end{corollary}
\begin{proof}
The first statement follows by iterated application of Lemma \ref{rules}.
 The last statement follows since
tensor products and contractions of tensors are bounded in terms of
the bounds on the factors, and because
permutations are isometric. Note that $\abs{\nu}=1$ and
$\abs{S^{*_1}}=\abs{S}$ for an arbitrary tensor $S$. Moreover, restriction
to the boundary only decreases the norm of a tensor.
\end{proof}

\begin{corollary}
  If $M$ is a Riemannian manifold of (coordinate-free defined) bounded
  geometry, the same is true for its boundary.
\end{corollary}

\medskip
\noindent\textbf{Step 4: Proof of Theorem \ref{main}(b1).}

Suppose we are in the situation of \ref{main}(b1) with $p\in\boundary
M$ and normal collar coordinates
$(x_1,\dots,x_m)=\kappa_p^{-1}:U\to\reals^m$ around $p$. By our
convention $x_m$ is the boundary defining coordinate,
i.e.~$\partial_m|_{\boundary M}=\nu$.

First consider $\boundary M$ as a Riemannian $(m-1)$-dimensional
manifold of its own.
Corollary \ref{boundarybounds} shows that bounds on the covariant
derivatives $\nabla^j R$ and $\bar\nabla^j l$ give rise to bounds on
$\bar\nabla^j\bar R$ ($0\le j\le k$). As in Step 2 (applied to
$\boundary M$) construct the orthonormal frame $\{s_i\}_{1\le i\le
  m-1}$ of $T\boundary M$. Extend this to an orthonormal frame of
$TM|_{\boundary M}$ by setting $s_m:=\nu$. By parallel transport along
geodesics with initial speed $\nu$ we get a synchronous orthonormal
frame of $TM$ on the normal collar neighborhood. Define the dual frame
$\{\theta^i\}$, the Christoffel symbols $\Gamma^i_{jk}$ and
$\tilde\Gamma^i_{jk}$, the curvature coefficients $R^i_{jkl}$ and
$K^i_{jkl}$, and $a^i_j$ and $b^i_j$  in
exactly the same way as in Step 2. Note that 
 \eqref{R_K} and \eqref{tgamma_gamma} remain true.

Now we come to the differential equations which relate these
quantities. By construction, $\{s_i\}$ is parallel to $\partial_m$. This translates to 
\begin{equation}\label{ith}
c(\partial_m)\theta^i_j=0,\quad\text{i.e.}\quad \Gamma^i_{jm}=0\qquad\forall i,j
\end{equation}
($c(\partial_m)$ denotes contraction with $\partial_m$).
The Lie derivative along $\partial_m$ (denoted by $\partial_m$) acts
on differential forms via $\partial_m=c(\partial_m)d + d
c(\partial_m)$. Hence
\begin{equation*}\begin{split} \partial_m\theta^i_j &=
    c(\partial_m)d\theta^i_j \stackrel{\eqref{ith}}{=}
    c(\partial_m)(d\theta^i_j - \sum_{k=1}^m \theta^i_k\wedge
    \theta^k_j)\\
&= c(\partial_m)(\sum_{k,l}R^i_{jkl})dx_k\wedge dx_l = 2\sum_k R^i_{jmk} dx_k .
\end{split}\end{equation*}
On the other hand,
$\partial_m\theta^i_j =
\sum_k\partial_m(\Gamma^i_{jk})dx_k$. 
Hence, applying $D^\alpha$ yields
\begin{equation}\label{nu_gamma}
\partial_m(D^\alpha\Gamma^i_{jk}) = 2 D^\alpha R^i_{jmk}\quad\forall i,j,k .
\end{equation}

Additionally, we need an equation for $a^i_j$. We apply $\partial_m$
twice to the dual frame. Since $\partial_m=s_m$ we have
$c(\partial_m)\theta^i=\delta_{im}$ ($\delta$ the Kronecker symbol). Then
\begin{equation*}\partial_m \theta^i =c(\partial_m)d\theta^i +d
c(\partial_m)\theta^i = c(\partial_m)d\theta^i . \end{equation*}
The connection is torsion free. This means
\begin{equation}\begin{split} d\theta^i &=\sum_j \theta^i_j\wedge\theta^j\\
\implies & \partial_m(\theta^i) = \sum_j
c(\partial_m)
(\theta^i_j\wedge\theta^j) \stackrel{\eqref{ith}}{=}-\theta^i_m
\label{theta1} \\
\implies &\partial_m^2(\theta^i)=-\partial_m(\theta^i_m)=-2\sum_k
R^i_{mmk}dx_k.  
\end{split}\end{equation}
The left hand side can be computed in terms of $a^i_j$
\begin{equation}\label{theta2} \partial_m(\theta^i)
  =\sum_j\partial_m(a^i_j)dx_j
\quad\implies\quad
\partial_m^2(\theta^i)=\sum_j\partial_m^2(a^i_j)dx_j . \end{equation}
Equating coefficients, applying $D^\alpha$, and expressing $R^i_{jkl}$
in terms of $K^i_{jkl}$ yields
\begin{equation}\label{nu_a}
 \partial_m^2 D^\alpha a^i_j = -2\sum_{k,l}D^\alpha(K^i_{mkl}a^k_m
a^l_j) \qquad\forall i,j .
\end{equation}
For $\abs{\alpha}>0$ this is (for each point in the boundary) a
system of inhomogeneous linear ordinary
differential equations for $D^\alpha a^i_j$, with coefficients given
by partial derivatives of $K^i_{mkl}$ up to order $\abs{\alpha}$ and
of $a^i_j$ up to order $\abs{\alpha}-1$.

To make use of the differential equation \eqref{nu_gamma} and
\eqref{nu_a} we have to determine the initial values
(at $x_m=0$).

If $i,j<m$, $a^i_j|_{x_m=0}$ is given by application of Step 2 to
$\boundary M$, which is possible because of Corollary
\ref{boundarybounds}. In particular, we get bounds for these
functions. And by construction  $a^i_m=a^m_i=\delta_{im}$.
For the first derivative we have $\partial_m a^i_j =
-\Gamma^i_{mj}$ (this follows from \eqref{theta1} and \eqref{theta2}).

Next we compute
$\Gamma^i_{jk}$ on $\boundary M$. By definition
$\bar\nabla_{\partial_i}s_j=\nabla_{\partial_i}s_j-
l(\partial_i,s_j)\nu$ for $i,j<m$. If we define
$l_{ij}:=l(\partial_i,s_j)$ then for $j,k<m$
\begin{equation*} \Gamma^i_{jk}|_{x_m=0} =\begin{cases}
    \bar\Gamma^i_{jk}|_{x_m=0}; &i<m\\ 
                                l_{kj}; &i=m . \end{cases}
\end{equation*}
By \eqref{ith} $\Gamma^i_{jm}=0$ $\forall i,j$. To compute
                                $\Gamma^i_{mj}$ we use for $j<m$
\begin{gather*}
g(\nabla_{\partial_j}\partial_m,s_i) =
\partial_jg(\partial_m,s_i)-
g(\partial_m,\nabla_{\partial_j}s_i) \stackrel{\partial_m|_{\boundary
    M}=\nu}{=} - l(\partial_j,s_i)\quad\text{for }i<m\\
2g(\nabla_{\partial_j}\partial_m,\partial_m) =\partial_j
g(\partial_m,\partial_m) =0\\
\implies   \nabla_{\partial_j}\partial_m=
  - \sum_i l_{ji} s_i \qquad\text{on
  $\boundary M$}.
\end{gather*}
It follows for $i<m$
\begin{equation*} \Gamma^i_{mj}|_{x_m=0} = \begin{cases} - l_{ji}; & j<m \\ 0; &
j=m . \end{cases} \end{equation*}
The arguments given in the proof of Step 2 show that if bounds exist
on the covariant derivatives up to order $k$ of the second fundamental
form and of $R$ (hence by Corollary \ref{boundarybounds} also on $\bar
R$) then the initial values of \eqref{nu_gamma} and \eqref{nu_a},
namely $D^\alpha\Gamma^i_{jk}|_{x_m=0}$,
$D^\alpha a^i_j|_{x_m=0}$ and $\partial_m D^\alpha a^i_j|_{x_m=0}$ are
bounded, as long as $\alpha_m=0$ (if
$\alpha=(\alpha_1,\dots,\alpha_m)$). Later, we will by induction on
$\abs{\alpha}$ get bounds on the right had sides of \eqref{nu_gamma}
and \eqref{nu_a} (using the bounds on the initial values), 
giving in particular bounds on $\partial_m D^\alpha
\Gamma^i_{jk}|_{x_m=0}$ and $\partial_m^2 D^\alpha a^i_j|_{x_m=0}$
which are the initial values of the equation for
$D^\beta\Gamma^i_{jk}$ and $D^\beta a^i_j$ where
$\beta=(\alpha_1,\dots,\alpha_m+1)$. We will therefore, inductively,
get the required bounds 
\begin{equation}\label{displaycite}
  \abs{D^\gamma\Gamma^i_{jk}|_{x_m=0}}\le C, \quad \abs{
D^\gamma a^i_j|_{x_m=0}}\le C,\quad \abs{\partial_m D^\gamma
a^i_j|_{x_m=0}}\le C.
\end{equation}

We proceed with a bootstrap argument similar to the one in Step 2. We
have to find bounds on $g_{ij}$, $g^{ij}$ and their
derivatives. Because of \eqref{defa} it suffices to look at $a^i_j$
and $b^i_j$ and their derivatives. Bounds on $a^i_j$ and $b^i_j$ are
given by Lemma \ref{abounds}. As in the proof of Step 2, because
of Lemma \ref{inv_mat} we only have to control the derivatives of
$a^i_j$. We do this by induction on the order of these
derivatives. To carry out the induction step, we also have to control
the derivatives of $K^i_{jkl}$, $R^i_{jkl}$ and $\Gamma^i_{jk}$ (and
the initial values in \eqref{displaycite}).

To conclude the start of the induction, Lemma \ref{Rder} and the assumptions give bounds on $K^i_{jkl}$ and (since the length of $\partial_i$ is bounded by
Proposition \ref{rauch}) on $R^i_{jkl}$. Integrating Equation
\eqref{nu_gamma}, we find bounds for
$\Gamma^i_{jk}$ (depending also on the given width $R_1$ of the
normal boundary charts). 

Assume now by induction that we have bounds on $D^\alpha a_{ij}$, $D^\alpha
\Gamma^k_{ij}$, $D^\alpha K^i_{jkl}$, and $D^\alpha R^i_{jkl}$ for
$\abs{\alpha}\le r$. By Lemma \ref{Rder}, the assumed bound on
$\nabla^{r+1}R$ therefore gives
bounds on $\nabla^\gamma K^i_{jkl}$ ($\abs{\gamma}=r+1$). Because of these and the bounds
on $\nabla^\alpha K^i_{jkl}$, $\nabla^\alpha a^i_j$ ($\abs{\alpha}\le
r$) and on the initial values \eqref{displaycite}
we can apply
\cite[IV.4.2]{Hartmann(1964)} to \eqref{nu_a} to obtain bounds on
$\nabla^\gamma a^i_j$ ($\abs{\gamma}=r+1$). This in turn,
together with the relation \eqref{R_K} between $R^i_{jkl}$ and
$K^i_{jkl}$ yields bounds on $\nabla^\gamma R^i_{jkl}$
($\abs{\gamma}=r+1$). 

Now we integrate \eqref{nu_gamma} to get bounds on $\nabla^\gamma
\Gamma^i_{jk}$ ($\abs{\gamma}=r+1$) to finish the induction step and to conclude the
proof of Theorem \ref{main}.

The bounds we obtain, inductively, depend only on the
bounds we started with.

\section{Technical properties of manifolds of bounded geometry}
\label{sec:prop}

We use the notation of Definition \ref{basicnotation}.

\begin{lemma}\label{balls}
Let $(M^n,g)$ be a  Riemannian manifold with boundary and with bounded
geometry as in Definition \ref{coorddef}.
  We find $r_0>0$ such that  for all $r,s\le r_0$ the following holds:
\begin{enumerate}
\item  If
 $x,x'\in\boundary M$ and
$ Z(x',s,\frac{2 R_2}{3})\cap U_r(Z(x,\frac{R_1}{2},\frac{2
  R_2}{3}))\ne \emptyset$
 then $Z(x',s,\frac{2 R_2}{3})\subset Z(x,\frac{9 R_1}{10},\frac{2 R_2}{3})$.
\item 
We find $0<D_1(r)<D_2(r)$ for $r\ge 0$
 such that 
\begin{equation*}\begin{split}
 & D_1(r)\le \vol(B(x,r))\le D_2(r) \qquad\forall x\in
 M-N(\frac{R_2}{3}),\text{ if }r<R_3\\
 &D_1(r)\le \vol(Z(x',r,\frac{2R_2}{3}))\le D_2(r)\qquad\forall x'\in
 \boundary M,\quad\text{if }r<R_1 .
\end{split}\end{equation*}
\end{enumerate}
$D_1(r)$, $D_2(r)$ and $r_0$ can be chosen to depend only on the
bounded geometry constants.
\end{lemma}
\begin{proof}
Bounded geometry implies the existence of $C_1,C_2>0$ so that in normal coordinates
\begin{equation*} \norm{(g^{ij})_{i,j}}<C_1,
\quad\norm{(g_{ij})_{i,j}}<C_1,\quad\text{and } C_2\le \sqrt{\abs{\det
(g_{ij})}}\le C_1 .
\end{equation*}
Observe that
$d(Z(x,\frac{R_1}{2},\frac{2R_2}{3}), M-Z(x,\frac{3R_1}{4},\frac{9R_2}{10}))$
is bounded independent of $x\in\boundary M$, using the bounds on the metric
tensor. With all sets and distances in $\boundary M$, $d(B(x,\frac{3R_1}{4}), \boundary
M- B(x,\frac{9R_1}{10}))\le R_1/10$. Choose $r_0$ smaller than half the minimum of these two
bounds. If $r,s<r_0$ and $Z(x',s,\frac{2R_2}{3})\cap
U_r(Z(x,\frac{R_1}{2},\frac{R_2}{2}))\ne\emptyset$ for
$x,x'\in\boundary M$ then
$Z(x',s,\frac{2R_2}{3})\cap
Z(x,\frac{3R_1}{4},\frac{9R_2}{10})\ne\emptyset$ which in turn implies
$Z(x',s,\frac{2R_2}{3})\subset Z(x,\frac{9R_1}{10},\frac{2 R_2}{3})$ by the choice of
$r_0$. This proves the first assertion.

The assertion about the volume bounds follows immediately from the
upper and lower
bounds of $\sqrt{\abs{\det(g_{ij})}}$. 

We can choose all constants to
depend only on the bounded geometry constants.
\end{proof}

The following is important to do analysis on
manifolds of bounded geometry. The corresponding result for empty
boundary is due to Shubin \cite[A1.2
  and A1.3]{Shubin(1992b)}.

\begin{proposition} (Partition of unity)\label{partition}\\
Let $M$ be a manifold with boundary and of bounded
geometry as in Definition \ref{coorddef}. There are $r_m>0$ and, for $0<r<r_m$ constants $C_K>0$
($K\in\naturals$), $M_f\in\naturals$, all depending only on the
bounded geometry constants (and $r$) such that
a covering of $M$ exists by sets $\{U(x_i,r)\}_{i\in
  I\subset\integers}$ which has the following properties:
\begin{enumerate}
\item\label{cover} $x_i\in\boundary M$ for $i\ge 0$ and $U(x_i,r)=
  Z(x_i,r,\frac{2R_2}{3})$;\\
  $x_i\in
  M-N(\frac{R_2}{2})$ for $i<0$ and $U(x_i,r)=B(x_i,r)$.
\item\label{locfin}  If $ s<r_m$ and
  $x\in M$ then $B(x,s)\cap U(x_i,r)\ne\emptyset$ for at most $M_f$
  of the $x_i$.
\item $\{U(x_i,r/2)\}_{i\in I}$ is a covering of $M$, too.
\end{enumerate}
Denote with $\kappa_i:B(0,r)\to U(x_i,r)$ ($i<0$) and
$\kappa_i:B(0,r)\times [0, \frac{2 R_2}{3})\to U(x_i,r)$ for $i\ge 0$ the
corresponding normal charts.

To this covering, a subordinate  partition of unity $\{\varphi_i\}$
exists such that
\begin{equation*} |D^\alpha\varphi_i|\le C_K\qquad\forall i\in\integers\quad\forall\abs{\alpha}\le K\qquad\text{(in normal coordinates).}\end{equation*}
\end{proposition}

\begin{proof} 
Set $r_m:=\min\{R_1/2,R_2/12, R_3, r_0/2\}$, where $r_0$ is given by Lemma
\ref{balls}.
Let $0<r<r_m$. First choose a maximal set of points $\{x_i\in\boundary
M;\; i=0,1,2,\dots\}$ such that all $(B(x_i,r/4)\in\boundary M)$ are
disjoint. Next, choose a maximal set of points $\{x_i\in M-
N(\frac{R_2}{2});\; i=-1,-2,\dots\}$ such that all $B(x_i,r/4)$
($i<0$) are disjoint. Note that the set $I$ of $i$ obtained this way
may be a proper subset of $\integers$. For $0<s\le r_0$ set $U(x_i,s):= B(x_i,s)$
if $i<0$
and $U(x_i,s):=Z(x_i,s,\frac{2R_2}{3})$ if $i\ge 0$. Then
\begin{equation*}
  \bigcup_{i<0}U(x_{i},r/2)=\bigcup_{i<0}B(x_i,r/2)\}\quad\text{covers
    $M-N(\frac{R_2}{2})$}.
 \end{equation*}
 This is true because else we find $z\in M-N(\frac{R_2}{2})$ which has distance
 $\ge r/2$ to all of the $x_{i}$. Then $B(z,r/4)\cap
 B(x_i,r/4)=\emptyset$ $\forall i<0$, violating the maximality of
 $\{x_{i}\}_{i<0}$. Similarly, $\{B(x_{i},r/2)\subset\boundary M\}_{i\ge 0}$ covers
 $\boundary M$ $\implies$ $\{U(x_{i},r/2)\}_{i\ge 0}$ covers $N(\frac{2R_2}{3})$.\\
Now we have to show that the covering $\{U(x_i,r)\}_{i\in I}$ has
 Property \ref{locfin}. So fix $0<s<r_m$ and $x\in M$.
 \begin{itemize}
 \item If $x\in N(\frac{R_2}{3})$ and $i<0$ then $B(x,s)\cap U(x_i,r)=\emptyset$
   since $d(N(\frac{R_2}{3}),M-N(\frac{R_2}{2})) =
   \frac{R_2}{6}>r+s$. 
 \item If $x\in M-N(\frac{R_2}{3})$ then the number $N_1$ of $x_i$
  ($i<0$) with $U(x_i,r)\cap B(x,s)\ne\emptyset$ is by Lemma
   \ref{balls} bounded by
   \begin{equation*}
     N_1\le \frac{\vol(B(x,s+r))}{\inf_{x_i\in
     M-N(\frac{R_2}{2})}\vol(B(x_i,r/4))} \le \frac{D_2(2r_m)}{D_1(r/4)}
\end{equation*}
since for such $x_i$ we have $B(x_i,r/4)\subset B(x,s+r)$ and all of
these are disjoint.
\item If $x\in M-N(R_2)$ then $B(x,s)\cap U(x_i,r)=\emptyset$ for
  $i\ge 0$ since $d(N(\frac{2R_2}{3}),M-N(R_2)) = \frac{R_2}{3}>s$.
\item If $x\in N(R_2)$ then the number $N_2$ of $x_i$ ($i\ge 0$) with
  $B(x,s)\cap U(x_i,r)\ne\emptyset$ is bounded by
  \begin{equation*}
    N_2\le \frac{\sup_{x_i\in\boundary
    M}\vol(Z(x_i,\frac{9R_1}{10},\frac{2 R_2}{3}))}{\inf_{x_i\in\boundary
    M}\vol(Z(x_i,\frac{r}{4},\frac{2R_2}{3}))} \le \frac{D_2(9R_1/10)}{D_1(r/4)}
\end{equation*}
since if there is one such $i_0$ then for all other such $i$ by Lemma \ref{balls}
$Z(x_i,\frac{r}{4},\frac{2R_2}{3})\subset Z(x_{i_0},9R_1/10,\frac{2 R_2}{3})$, and all
the $Z(x_i,\frac{r}{4},\frac{2R_2}{3})$ are disjoint.

 It follows in
all cases
\begin{equation*}
  M_f(r) \le \frac{D_2(9R_1/10)+D_2(2r_0)}{D_1(r/4)}< \infty.
\end{equation*}
 \end{itemize}

It remains to construct the subordinate partition of unity. Choose a
smooth cutoff function $\varphi:\reals^m\to [0,1]$ with $\varphi(x)=1$ if
$\abs{x}\le r/2$ and $\varphi(x)=0$ if $\abs{x}\ge r$. Denote the
restriction to $\reals^{m-1}$ also with $\varphi$. Choose smooth
$\psi:\reals\to [0,1]$ with $\psi(x)=0$ if $x\ge 2R_2/3$ and $\psi(x)=1$ if
$x\le R_2/2$. Via the normal coordinates this yields cutoff functions
$f_i$ on $U(x_i,r)$ with $f_{i}\circ\kappa_i(y',t)=\varphi(y')\psi(t)$ if
$i\ge 0$ and $f_{i}\circ\kappa_i=\varphi$ if $i<0$. Therefore, if $\kappa$
is any normal chart,
$f_{i}\circ\kappa=\varphi\circ(\kappa_{i}^{-1}\circ\kappa)$ ($i<0$) and
$f_{i}\circ\kappa=(\varphi\cdot\psi)\circ(\kappa^{-1}_i\circ\kappa)$ ($i\ge
0$). The chain rule shows that the  bounds on derivatives up to order
$K$ of the coordinate changes (Proposition \ref{change_var}) yield
bounds on the partial derivatives up to order $K$ of $f_i$ in normal
coordinates. To construct the partition of unity, set $F=\sum_{i\in I}
f_i$ (at each point there are at most $M_f$ non-zero summands). 
\begin{equation*}
\text{Since}\qquad  M-N(\frac{R_2}{2})\subset \bigcup_{i<0}
 U(x_{i},r/2)\quad\text{and}\quad N(\frac{R_2}{2})\subset\bigcup_{i\ge
 0} U(x_{i},r/2), 
\end{equation*}
 for each $z\in M$ at least one of $f_i(z)=1$ $\implies$ $F\ge 1 $. Define 
\begin{equation*}\varphi_i:=f_i/F.\end{equation*}
Obviously, $\{\varphi_i\}_{i\in I}$ is a smooth partition of unity
subordinate to our covering. Pick one $\varphi_i$ and one normal chart $\kappa$.  For partial derivatives up to order $K$ in normal coordinates observe
\begin{equation*}\begin{split} \abs{D^\alpha(\varphi_i\circ\kappa)}
  &=\abs{D^\alpha\frac{f_i\circ\kappa}{F\circ\kappa}}
   =\frac{\abs{P_\alpha(D^\beta(f_i\circ\kappa),D^\gamma(F\circ\kappa);\;\abs{\beta},\abs{\gamma}\le\abs{\alpha})}}{\abs{F\circ\kappa}^{2^{\abs{\alpha}}}}\\
&\stackrel{\abs{F}\ge 1}{\le}
\abs{P_\alpha(D^\beta(f_i\circ\kappa),D^\gamma(F\circ\kappa))} . 
\end{split}\end{equation*}
$P_\alpha$ is a polynomial entirely determined by $\alpha$. At every
point $x\in M$, $D^\gamma(F\circ\kappa)|_x$ is the sum of at most $M_f$
summands of the type $D^\gamma(f_s\circ\kappa)|_x$. Therefore, we have
bounds for all the entries of $P_\alpha$. This yields a bound $C_K$,
depending only on the bounded geometry constants, for
$\abs{D^\alpha(\varphi_i\circ\kappa)}$ if $\abs{\alpha}\le K$.  
\end{proof}

\medskip\noindent
\textbf{Changes of normal coordinates}

\begin{proposition}\label{coord_change}\label{change_var}
Suppose $M$ is a  Riemannian manifold with boundary and of bounded
geometry. More precisely, suppose\ $C>0$ is a bound for partial
derivatives up to order $k+1$ of $g^{ij}$ and $g_{ij}$ in normal
coordinates. Then $D>0$ exists, depending only on $C$ so that, if $\kappa_1:U_1\subset\reals^m\to
M$ and $\kappa_2:U_2\subset\reals^m\to M$ are  normal charts as in \ref{main}, the
following holds for $f:=\kappa_1^{-1}\circ\kappa_2: U_0\subset\reals^m\to
\reals^m$ ($U_0$ the domain of definition of the composition):
\begin{equation*} \abs{D^\alpha f}\le D\quad\forall\abs{\alpha}\le k .\end{equation*}
\end{proposition}

Since the maps $\kappa_i$ are solutions of certain ordinary
differential equation, namely the equation for geodesics, we first
recall a result about differential equations.

\begin{lemma}\label{sol_i}
Let $x'(t)=F(t,x(t))$ be a system of ordinary differential equations
($t\in\reals$, $x(t)\in\reals^n$), $F\in C^\infty(\reals\times\reals^n,\reals^n)$. Suppose
$\varphi(t,x)$ is the flow of this equation. We find a universal
expression $Expr_\alpha$,
only depending on $\alpha$ such
that  for all $ t\ge 0$ where $\phi(t,x_0)$ is defined
\begin{equation}\label{ie1} 
\abs{D_x^\alpha \phi(t,x_0)} \le Expr_\alpha\left(
\sup_{0\le\tau\le t}\{\abs{D_x^\beta F(\tau,\phi(\tau,x_0)}\}
|\;\beta\le\alpha, t\right) . 
\end{equation}
\end{lemma}
\begin{proof}
The theory of ordinary differential equations\ \cite[V.3.1.]{Hartmann(1964)} tells us   that we have the linear differential equation
\begin{equation*} \alpha'(t)=\frac{\partial F}{\partial
  x}(t,\varphi(t,x))\cdot\alpha(t);\quad
\alpha(0)=e_k=(0,\dots,1,\dots,0) \end{equation*}
for $\partial_k\varphi(t,x)$. For linear differential equations \cite[IV.4.2]{Hartmann(1964)}
gives inequalities  which
directly imply \eqref{ie1} if $\abs{\alpha}=1$.
  
Inductively one shows that
higher derivatives fulfill the linear differential equation
\begin{multline}
  (D_x^\alpha\varphi)'(t,x)=(D_xF)(t,\varphi(t,x))\cdot
  D_x^\alpha\varphi(t,x) +\\ P_\alpha(D_x^\gamma\varphi, (D_x^\beta
  F)(t,\varphi(t,x));\; \gamma<\alpha,\; \beta\le\alpha)
\end{multline}
with $D_x^\alpha\varphi(0,x)= 0$ if $\abs{\alpha}>1$. 
Here $P_\alpha$ is a polynomial matrix which depends only on
$\alpha$. By induction and using \cite[IV.4.2]{Hartmann(1964)} again, the
proposition follows.
\end{proof}

Reduction of order implies:
\begin{corollary}\label{sol_ie}
A statement corresponding to Lemma \ref{sol_i} holds for
ordinary differential equations of order $k$.
\end{corollary}

\begin{lemma}
Let $p\in U\subset\reals^m$, $g_{ij}$ a Riemannian metric on $U$ and
$\exp_p:B(r,0)\to U$ the exponential map at $p$ (we identify $T_pU$
with $\reals^m$ via an orthonormal frame). If $r$ is sufficiently small
then $\exp_p$ is a diffeomorphism onto some open subset of $U$, and
the derivatives up to order $k$ of $\exp_p$ and its inverse are
bounded in terms of $g_{ij}$, $g^{ij}$, their derivatives up to order $k+1$, and $r$.
\end{lemma}
\begin{proof}
We have $\exp_p(x)=\varphi(x,p,1)$, where $\varphi$ with
$\varphi(x,q,0)=q$, $\varphi'(x,q,0)=x$ is the flow of the differential equation for geodesics. 
Corollary \ref{sol_ie} applies to this equation $x''=F(x)$, and
$F(x)=-\sum_{i,j}\Gamma_{ij}(x)x_i'x_j'$ is given by $g_{ij}$ and its
first order derivatives.

For the inverse, by Lemma \ref{inv_mat} it suffices to study its first
order derivatives. Bounds on these follow from Proposition \ref{rauch}.
\end{proof}

\begin{lemma}
Suppose $ U',V\subset\reals^{m-1}$, $\kappa:U'\times [0,r_C)\to
V\times [0,r_C)$ is a normal boundary chart centered at $p\in V$ on the Riemannian manifold
$V\times [0,r_C)$ with metric $g_{ij}$
($g_{im}=\delta_{im}=g_{mi}$). Then the derivatives of $\kappa$ and
its inverse up to order $k$ are bounded in terms of $g_{ij}$, $g^{ij}
$ and their derivatives up to order $k+1$.
\end{lemma}
\begin{proof}
$\kappa(q,s)=\varphi_1(s\cdot \partial_m,\varphi_2(q,p,1),0,1)$, where
$\varphi_1$ is the flow of the differential equation for the geodesics
in $V\times [0,r_C)$ ($\varphi_1(v,p,\tau,0)=(p,\tau)$,
$\varphi_1'(v,(p,\tau),0)=v$), and $\varphi_2$ is the flow of the
differential equation for geodesics on $V$. Hence $\kappa$ is the
composition of two flows to which Corollary \ref{sol_ie}, and then Lemma
\ref{inv_mat} and Proposition \ref{rauch} applies exactly as in 
the previous lemma.
\end{proof}

We prove Proposition \ref{coord_change} using these Lemmas as
follows: 
By Theorem \ref{main} we have  bounds for $g_{ij}$ and their derivatives up to order $k+1$ in normal coordinates. Write
\begin{equation*} \kappa_1^{-1}\circ \kappa_2 =
  (\kappa_1^{-1}\circ\kappa_0)\circ
  (\kappa_0^{-1}\circ\kappa_2). \end{equation*}
with $\kappa_0$ either being an exponential map or a normal boundary
  map with suitable range, respectively. If we use $\kappa_1$ or
  $\kappa_2$ to pull back the given Riemannian metric to the domain of
  the charts, $\kappa_0^{-1}\circ\kappa_2$ and
  $\kappa_1^{-1}\circ\kappa_0$ each fulfill exactly the assumptions of
  one of the two lemmas. The conclusion of these lemmas is then true
  for there composition, as well, and the Proposition follows.

\begin{corollary}
  To check condition (B1) in Definition \ref{coorddef} it suffices to
  do this for an atlas of such charts.
\end{corollary}



\end{document}